\documentclass{article}
\usepackage{amsmath,amsfonts,amsthm}
\usepackage{amstext}
\usepackage{amsgen}
\usepackage{amsbsy}
\usepackage{amsopn}
\usepackage{amssymb}
\usepackage{graphicx}
\usepackage{epsfig}

\theoremstyle{definition}

\theoremstyle{remark}

\numberwithin{equation}{section}

\begin{document}

\title{Conjugate Gradient Algorithm for Solving a Optimal Multiply Control
Problem on \\
a System of Partial Differential Equations}
\author{Carlos Barr\'{o}n Romero \\
cbarron@correo.azc.uam.mx\\
UAM-Azcapotzalco \\
Department of Basic Sciences}
\date{}
\maketitle

\begin{abstract}
I development a Conjugate Gradient Method for solving a partial
differential system with multiply controls. Some numerical results
are depicted. Also, I present an explication of why the control
over a partial differential equations system is necessary.

\textbf{Keywords}: Optimal Control over Partial Differential Equations;
Process Engineering Methods.
\end{abstract}





\section{Introduction}

Given the partial differential system:

\begin{center}
\begin{equation}
\left\{
\begin{tabular}{ll}
$\frac{\partial y}{\partial t}-\mu \frac{\partial ^{2}y}{\partial x^{2}}%
+\epsilon \frac{\partial y}{\partial x}-y=0$ & $\text{in }Q=\left(
0,L\right) \times \left( 0,T\right) $ \\
$y\left( x,0\right) =y_{0},$ & $t=0,$ \\
$-\mu \frac{\partial y\left( 0,t\right) }{\partial x}=0,$ & $x=0,$ \\
$\mu \frac{\partial y\left( L,t\right) }{\partial x}=0,$ & $x=L.$%
\end{tabular}
\right.  \tag{S}  \label{S}
\end{equation}
\end{center}

A conjugate gradient algorithm with several control on $\left[ 0,L\right] $
is developed for ~\ref{S}, which is similar to the Burgers' equation.

\begin{figure}
\centerline{ \psfig{figure=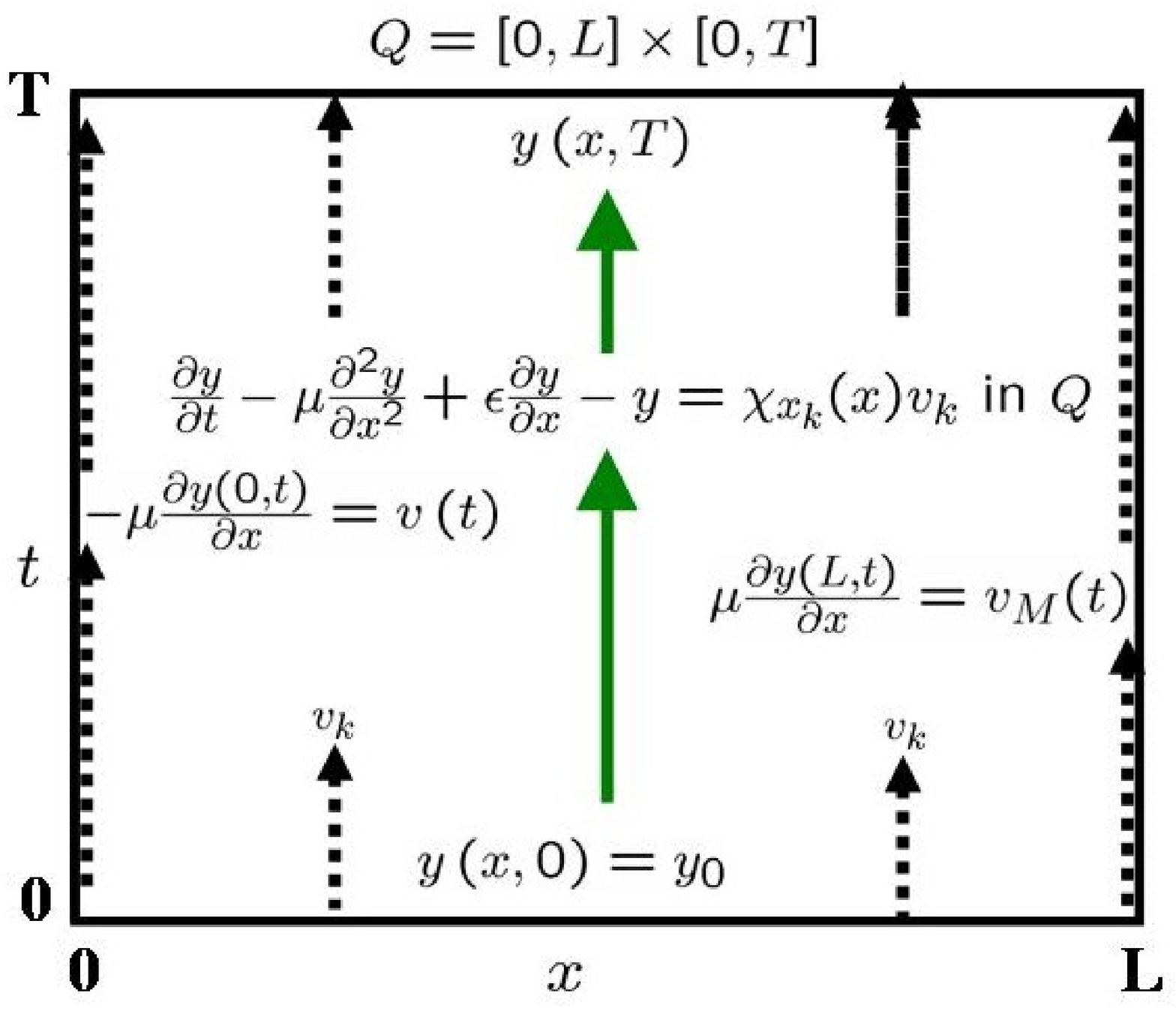,width=3.5in}}~ \caption{System
(\ref{SE}).} \label{fig:SE}
\end{figure}

\section{Several Control for ~\ref{S}}

With an appropriate functions $v \in \mathcal{V}$, $v=(v_0,v_1,\ldots,v_M)$
and $\mathcal{V}$ appropriate Hilbert Space, the system can be controlled on
$x_{k}=L\frac{k}{M},$ $k=0,\ldots ,M$\ (see figure \ref{fig:SE}).

\begin{center}
\begin{equation}
\left\{
\begin{tabular}{ll}
$\frac{\partial y}{\partial t}-\mu \frac{\partial ^{2}y}{\partial x^{2}}%
+\epsilon \frac{\partial y}{\partial x}-y$ $=\chi_{x_k}(x)v_k$ & $\text{in }%
Q=\left( 0,L\right) \times \left( 0,T\right) ,k=1,\ldots ,M-1$ \\
$y\left( x,0\right) =y_{0},$ & $t=0,$ \\
$-\mu \frac{\partial y\left( 0,t\right) }{\partial x}=v_{0},$ & $x=0,$ \\
$\mu \frac{\partial y\left( L,t\right) }{\partial x}=v_{M},$ & $x=L.$%
\end{tabular}
\right.  \tag{SE}  \label{SE}
\end{equation}
\end{center}

In this case, the corresponding variational control problem is
\begin{equation}
\left\{
\begin{tabular}{l}
Find $u^{\ast }\in \mathcal{V}$, \\
$J\left( u^{\ast }\right) \leq J\left( v\right) $, $\forall v\in \mathcal{V}$%
\end{tabular}%
\ \right.  \tag{CP}  \label{CP}
\end{equation}

where

\begin{equation*}
J\left( v\right) =\frac{k_{0}}{2}\sum_{k=0}^{M}\int_{0}^{T}v_{k}^{2}\text{d}%
t+\frac{k_{1}}{2}\iint_{Q}y^{2}\text{d}x\text{d}t+\frac{k_{2}}{2}%
\int_{0}^{L}y\left( x,T\right) ^{2}\text{d}x
\end{equation*}

where $v=(v_{0},v_{1},\ldots ,v_{M}),$ and $y$ is the solution of (\ref{SE})
for each $v$ (see figure \ref{fig:SE}).

The equivalent form as an optimization problem is:

\begin{equation*}
\min_{v\in \mathcal{U}}J\left( v\right) =\frac{k_{0}}{2}\sum_{k=0}^{M}%
\int_{0}^{T}v_{k}^{2}\text{d}t+\frac{k_{1}}{2}\iint_{Q}y^{2}\text{d}x\text{d}%
t+\frac{k_{2}}{2}\int_{0}^{L}y\left( x,T\right) ^{2}\text{d}x,
\end{equation*}

where $y$ is the solution of (\ref{SE}) for $v.$

In this case, the objective of the optimization problem is given a
perturbation function $y_0$ at $t=0$ get back to the steady state
to  $\mathbf{0}$. Also, the controls must reduce the cost or
weight of control variable $v$, keep low the cost of the evolution
of the system $y\left( x,t\right).$

\section{The continuous case}

The continuous case is computing by a perturbation of (\ref{CP}) and (\ref%
{SE}) and using the optimal (necessary and sufficient) condition $\delta
J\left( v\right) =0.$

\begin{eqnarray*}
\delta J\left( v\right)  &=&k_{0}\sum_{k=0}^{M}\int_{0}^{T}v_{k}\delta v_{k}%
\text{d}t+k_{1}\iint_{Q}y\delta y\text{d}x\text{d}t \\
&+&k_{2}\int_{0}^{L}y\left( x,T\right) \delta y\left( x,T\right) \text{d}x.
\end{eqnarray*}

The perturbation system of the equation (\ref{SE}) is

\begin{equation}
\left\{
\begin{tabular}{ll}
$\frac{\partial \delta y}{\partial t}-\mu \frac{\partial ^{2}\delta y}{%
\partial x^{2}}+\epsilon \frac{\partial \delta y}{\partial x}-\delta y=
\chi_{x_k}(x) \delta v_{k}$ & in $Q=\left( 0,L\right) \times \left(
0,T\right) ,k=1,\ldots ,M-1$ \\
$\delta y\left( x,0\right) =0,$ & $t=0,$ \\
$-\mu \frac{\partial \delta y\left( 0,t\right) }{\partial x}=\delta v_{0},$
& $x=0,$ \\
$\mu \frac{\partial \delta y\left( L,t\right) }{\partial x}=\delta v_{M},$ &
$x=L.$%
\end{tabular}
\right.  \tag{$\delta $SE}  \label{delSE}
\end{equation}

Let $p\left( x,t\right) $ a sufficiently smooth function \ that allow to
integrate ($\delta $SE) in $Q$
\begin{eqnarray*}
0 &=&\iint_{Q}p\left( \frac{\partial \delta y}{\partial t}-\mu \frac{%
\partial ^{2}\delta y}{\partial x^{2}}+\epsilon \frac{\partial \delta y}{%
\partial x}-\delta y-\chi _{x_{i}}\delta v\right) \text{d}x\text{d}t \\
&=&\iint_{Q}p\frac{\partial \delta y}{\partial t}\text{d}x\text{d}t-\mu
\iint_{Q}p\frac{\partial ^{2}\delta y}{\partial x^{2}}\text{d}x\text{d}%
t+\epsilon \iint_{Q}p\frac{\partial \delta y}{\partial x}\text{d}x\text{d}t
\\
&-& \iint_{Q}p\delta y\text{d}x\text{d}t-\iint_{Q}p\chi _{x_{i}}\delta v%
\text{d}x\text{d}t.
\end{eqnarray*}

The integration of ($\delta $SE) is achieved by the formula of integration
by parts:

\begin{equation*}
\int_{a}^{b}v\text{d}u = vu|_{a}^{b}-\int_{a}^{b}u\text{d}v.
\end{equation*}

Therefore

\begin{eqnarray}
\iint_{Q}p\frac{\partial \delta y}{\partial t}\text{d}x\text{d}t
&=&\int_{0}^{L}\left[ \int_{0}^{T}p\frac{\partial \delta y}{\partial t}\text{%
d}t\right] \text{d}x~  \label{1} \\
v &=&p,\ \text{d}u=\frac{\partial \delta y}{\partial t}\text{d}t  \notag \\
&=&\int_{0}^{L}\left[ p\left( x,T\right) \delta y\left( x,T\right) \right]
_{0}^{T}\text{d}x-\iint_{Q}\frac{\partial p}{\partial t}\delta y\text{d}x%
\text{d}t  \notag \\
&=&\int_{0}^{L}p\left( x,T\right) \delta y\left( x,T\right) \text{d}%
x-\int_{0}^{L}p\left( x,0\right) \delta y\left( x,0\right) \text{d}x  \notag
\\
&-& \iint_{Q}\frac{\partial p}{\partial t}\delta y\text{d}x\text{d}t  \notag
\\
&&\left( \delta y\left( x,0\right) =0\right)  \notag \\
&=&\int_{0}^{L}p\left( x,T\right) \delta y\left( x,T\right) \text{d}%
x+\iint_{Q}\left( -\frac{\partial p}{\partial t}\right) \delta y\text{d}x%
\text{d}t  \notag
\end{eqnarray}

\begin{equation}
-\mu \iint_{Q}p\frac{\partial ^{2}\delta y}{\partial x^{2}}\text{d}x\text{d}
t = -\mu \int_{0}^{T}\left[ \int_{0}^{L}p\frac{\partial ^{2}\delta y}{
\partial x^{2}}\text{d}x\right] \text{d}t~  \label{2}
\end{equation}

\begin{eqnarray}
v &=&p,\ \text{d}u=\frac{\partial ^{2}\delta y}{\partial x^{2}}\text{d}x
\notag \\
&=&-\mu \int_{0}^{T}\left[ p\left( x,t\right) \frac{\partial \delta y\left(
x,t\right) }{\partial x}\right] _{0}^{L}\text{d}t+\mu \iint_{Q}\left[ \frac{%
\partial p}{\partial x}\frac{\partial \delta y}{\partial x}\text{d}x\right]
\text{d}t  \notag \\
v &=&\frac{\partial p}{\partial x},\ \text{d}u=\frac{\partial \delta y}{%
\partial x}\text{d}x  \notag \\
&=&\int_{0}^{T}p\left( L,t\right) \left( -\mu \frac{\partial \delta y\left(
L,t\right) }{\partial x}\right) \text{d}t-\int_{0}^{T}p\left( 0,t\right)
\left( -\mu \frac{\partial \delta y\left( 0,t\right) }{\partial x}\right)
\text{d}t  \notag \\
&&+\mu \int_{0}^{T}\left[ \frac{\partial p\left( x,t\right) }{\partial x}%
\delta y\left( x,t\right) \right] _{0}^{L}\text{d}t-\mu \iint_{Q}\frac{%
\partial ^{2}p}{\partial x^{2}}\delta y\text{d}x\text{d}t.  \notag
\end{eqnarray}

\begin{eqnarray*}
&&\left( \mu \frac{\partial \delta y\left( L,t\right) }{\partial x}=\delta
v_{M}\left( t\right) ,-\mu \frac{\partial \delta y\left( 0,t\right) }{%
\partial x}=\delta v_{0}\left( t\right) \right) \\
&=&-\int_{0}^{T}p\left( L,t\right) \delta v_{M}\left( t\right) \text{d}%
t-\int_{0}^{T}p\left( 0,t\right) \delta v_{0}\left( t\right) \text{d}t \\
&&+\mu \int_{0}^{T}\left[ \frac{\partial p\left( x,t\right) }{\partial x}%
\delta y\left( x,t\right) \right] _{0}^{L}\text{d}t - \mu \iint_{Q}\frac{%
\partial ^{2}p}{\partial x^{2}}\delta y\text{d}x\text{d}t \\
&=&\int_{0}^{T}p\left( L,t\right) \left( -\delta v_{M}\left( t\right)
\right) \text{d}t-\int_{0}^{T}p\left( 0,t\right) \delta v_{0}\left( t\right)
\text{d}t \\
&& +\mu \int_{0}^{T}\frac{\partial p\left( L,t\right) }{\partial x}\delta
y\left( L,t\right) \text{d}t -\mu \int_{0}^{T}\frac{\partial p\left(
0,t\right) }{\partial x}\delta y\left( 0,t\right) \text{d}t \\
&&-\mu \iint_{Q}\frac{\partial ^{2}p}{\partial x^{2}}\delta y\text{d}x\text{d%
}t \\
&=&\int_{0}^{T}\left( -p\left( L,t\right) \right) \left( \delta v_{M}\left(
t\right) \right) \text{d}t+\int_{0}^{T}\left( -p\left( 0,t\right) \right)
\delta v_{0}\left( t\right) \text{d}t \\
&& +\int_{0}^{T}\mu \frac{\partial p\left( L,t\right) }{\partial x}\delta
y\left( L,t\right) \text{d}t +\int_{0}^{T}\left( -\mu \frac{\partial p\left(
0,t\right) }{\partial x}\right) \delta y\left( 0,t\right) \text{d}t \\
&& -\mu \iint_{Q}\frac{\partial ^{2}p}{\partial x^{2}}\delta y\text{d}x\text{%
d}t \\
&=&\int_{0}^{T}\left( -p\left( L,t\right) \right) \left( \delta v_{M}\left(
t\right) \right) \text{d}t+\int_{0}^{T}\left( -p\left( 0,t\right) \right)
\delta v_{0}\left( t\right) \text{d}t \\
&& +\int_{0}^{T}\mu \frac{\partial p\left( L,t\right) }{\partial x}\delta
y\left( L,t\right) \text{d}t + \int_{0}^{T}\left( -\mu \frac{\partial
p\left( 0,t\right) }{\partial x}\right) \delta y\left( 0,t\right) \text{d}t
\\
&& +\iint_{Q}\left( -\mu \frac{\partial ^{2}p}{\partial x^{2}}\right) \delta
y\text{d}x\text{d}t
\end{eqnarray*}

\begin{eqnarray}
\epsilon \iint_{Q}p\frac{\partial \delta y}{\partial x}\text{d}x\text{d}t
&=&\epsilon \int_{0}^{T}\left[ \int_{0}^{L}p\frac{\partial \delta y}{%
\partial x}\text{d}x\right] \text{d}t ~  \label{3}
\end{eqnarray}

\begin{eqnarray}
v &=&p,\ \text{d}u=\frac{\partial \delta y}{\partial x}\text{d}x  \notag \\
&=&\epsilon \int_{0}^{T}\left[ p\left( x,t\right) \delta y\left( x,t\right) %
\right] _{0}^{L}\text{d}t-\epsilon \iint_{Q}\frac{\partial p}{\partial x}%
\delta y\text{d}x\text{d}t  \notag \\
&=&\epsilon \int_{0}^{T}p\left( L,t\right) \delta y\left( L,t\right) \text{d}%
t-\epsilon \int_{0}^{T}p\left( 0,t\right) \delta y\left( 0,t\right) \text{d}t
\notag \\
&&- \epsilon \iint_{Q}\frac{\partial p}{\partial x}\delta y\text{d}x\text{d}t
\notag \\
&=&\int_{0}^{T}\epsilon p\left( L,t\right) \delta y\left( L,t\right) \text{d}%
t+\int_{0}^{T}\left( -\epsilon p\left( 0,t\right) \right) \delta y\left(
0,t\right) \text{d}t  \notag \\
&&- \epsilon \iint_{Q}\frac{\partial p}{\partial x}\delta y\text{d}x\text{d}t
\notag
\end{eqnarray}

\begin{equation}
-\iint_{Q}p\delta y\text{d}x\text{d}t=\iint_{Q}\left( -p\right) \delta y%
\text{d}x\text{d}t.~  \label{4}
\end{equation}

\begin{eqnarray}
-\iint_{Q}p\chi _{x_{i}}\delta v\text{d}x\text{d}t
&=&\sum_{i=1}^{M-1}\int_{0}^{T}\left( -p_{i}\right) \delta v_{i}\text{d}t~
\label{5} \\
\text{where }\chi _{x_{i}}p &=&p_{i}.~  \notag
\end{eqnarray}

\begin{eqnarray*}
0 &=&(\ref{1})+(\ref{2})+(\ref{3})+(\ref{4})+(\ref{5}) \\
&=&\int_{0}^{L}p\left( x,T\right) \delta y\left( x,T\right) \text{d}%
x+\iint_{Q}\left( -\frac{\partial p}{\partial t}\right) \delta y\text{d}x%
\text{d}t \\
&&+\int_{0}^{T}\left( -p\left( L,t\right) \right) \left( \delta v_{M}\left(
t\right) \right) \text{d}t+\int_{0}^{T}\left( -p\left( 0,t\right) \right)
\delta v_{0}\left( t\right) \text{d}t+\int_{0}^{T}\mu \frac{\partial p\left(
L,t\right) }{\partial x}\delta y\left( L,t\right) \text{d}t \\
&&+\int_{0}^{T}\left( -\mu \frac{\partial p\left( 0,t\right) }{\partial x}%
\right) \delta y\left( 0,t\right) \text{d}t+\iint_{Q}\left( -\mu \frac{%
\partial ^{2}p}{\partial x^{2}}\right) \delta y\text{d}x\text{d}t \\
&&+\int_{0}^{T}\epsilon p\left( L,t\right) \delta y\left( L,t\right) \text{d}%
t+\int_{0}^{T}\left( -\in p\left( 0,t\right) \right) \delta y\left(
0,t\right) \text{d}t-\epsilon \iint_{Q}\frac{\partial p}{\partial x}\delta y%
\text{d}x\text{d}t \\
&&+\iint_{Q}\left( -p\right) \delta y\text{d}x\text{d}t \\
&&+\sum_{i=1}^{M-1}\int_{0}^{T}\left( -p_{i}\right) \delta v_{i}\text{d}t \\
&=&\sum_{i=0}^{M}\int_{0}^{T}\left( -p_{i}\right) \delta v_{i}\text{d}t \\
&&+\iint_{Q}\left( -\frac{\partial p}{\partial t}-\mu \frac{\partial ^{2}p}{%
\partial x^{2}}-\epsilon \frac{\partial p}{\partial x}-p\right) \delta y%
\text{d}x\text{d}t \\
&&+\int_{0}^{L}p\left( x,T\right) \delta y\left( x,T\right) \text{d}x \\
&&+\int_{0}^{T}\left( \mu \frac{\partial p\left( L,t\right) }{\partial x}%
+\epsilon p\left( L,t\right) \right) \delta y\left( L,t\right) \text{d}%
t+\int_{0}^{T}\left( -\mu \frac{\partial p\left( 0,t\right) }{\partial x}%
-\epsilon p\left( 0,t\right) \right) \delta y\left( 0,t\right) \text{d}t
\end{eqnarray*}

Adjusting terms with
\begin{equation*}
\delta J\left( v\right) =k_{0}\sum_{k=0}^{M}\int_{0}^{T}v_{k}\delta v_{k}%
\text{d}t+k_{1}\iint_{Q}y\delta y\text{d}x\text{d}t+k_{2}\int_{0}^{L}y\left(
x,T\right) \delta y\left( x,T\right) \text{d}x,
\end{equation*}

the adjoint system is
\begin{equation}
\left\{
\begin{tabular}{ll}
$p\left( x,T\right) =k_{2}y\left( x,T\right) $, & $x\in \left[ 0,L\right] $
\\
$\mu \frac{\partial p}{\partial x}\left( L,t\right) +\epsilon p\left(
L,t\right) =0$, & $t\in \left[ 0,T\right] $ \\
$\mu \frac{\partial p}{\partial x}\left( 0,t\right) +\epsilon p\left(
0,t\right) =0$ & $t\in \left[ 0,T\right] $ \\
$\frac{\partial p}{\partial t}+\mu \frac{\partial ^{2}p}{\partial x^{2}}%
+\epsilon \frac{\partial p}{\partial x}+p=-k_{1}y,$ & in $Q$%
\end{tabular}%
\ \right.   \tag{$\delta $ASE}  \label{delASE}
\end{equation}%
also
\begin{equation*}
\nabla J\left( v\right) =k_{0}\sum_{k=0}^{M}\left( v_{k}-p_{k}\left(
x,t\right) \right) .
\end{equation*}

\begin{figure}
\centerline{ \psfig{figure=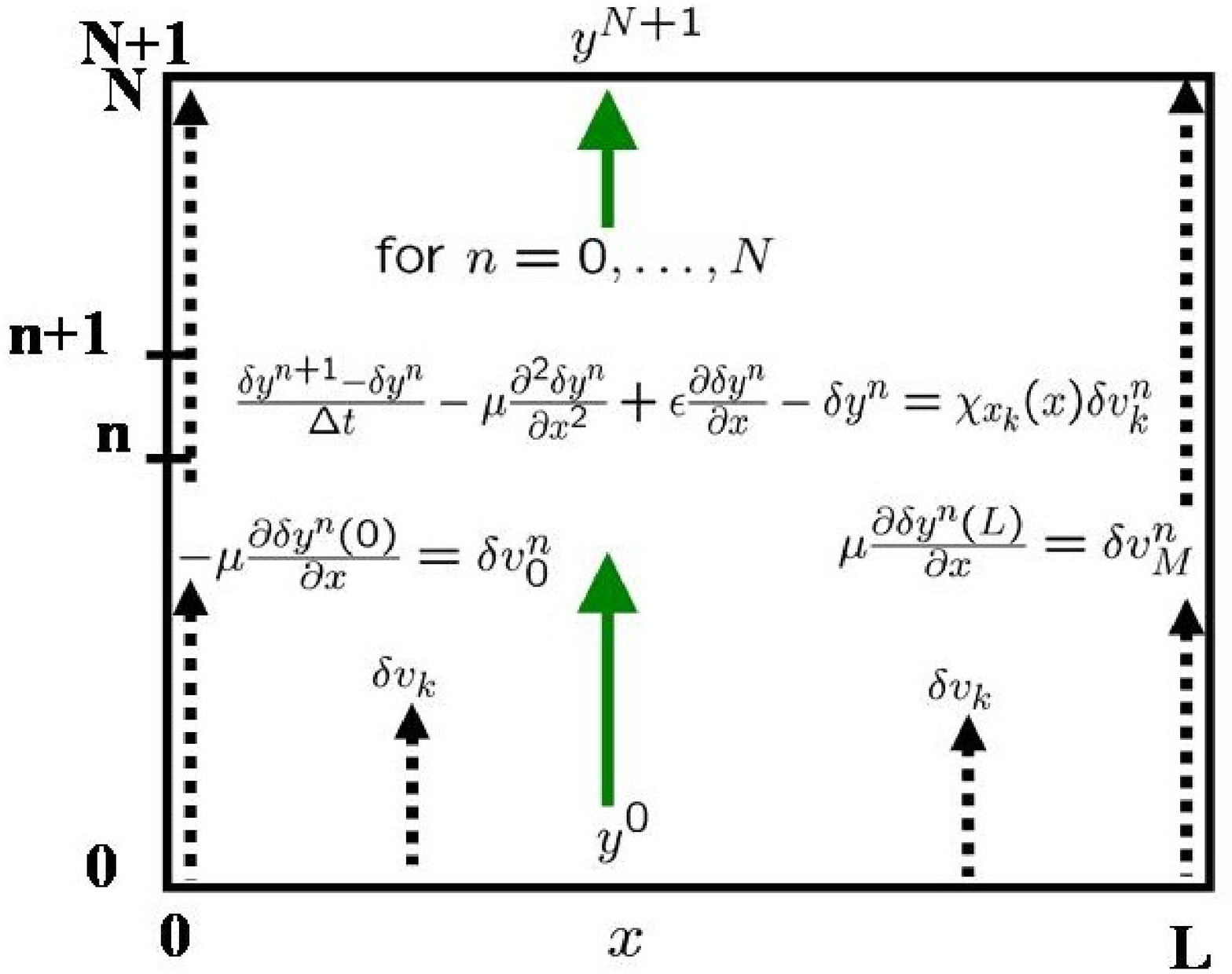,width=3.5in} } ~
\caption{Discretization on time of (\ref{delSEdelT}).}
\label{fig:delSE_DT}
\end{figure}

\section{Discretization on Time}

The discretization on time of $J^{^{{\Delta t}}}\left( v\right) $ is
\begin{equation*}
J^{^{{\Delta t}}}\left( v\right) =\frac{{\Delta t}}{2}\sum_{k=0}^{M}%
\sum_{n=0}^{N}\left\Vert v_k^{n}\right\Vert ^{2}+\frac{k_{1}{\Delta t}}{2}%
\sum_{n=0}^{N}\int_{0}^{L}\left\Vert y^{n}\right\Vert ^{2}\text{d}x+\frac{%
k_{2}}{2}\int_{0}^{L}\left\Vert y^{N+1}\left( x\right) \right\Vert ^{2}\text{%
d}x
\end{equation*}%
where $N>0$,and ${\Delta t}=\frac{T}{N}$.

Now, the forward discretization on time of (\ref{SE}) is

\begin{equation}
\left\{
\begin{tabular}{l}
$y^{0}=y_{0}.$ \\
for $n=0,\ldots ,N$ \\
$\frac{y^{n+1}-y^{n}}{{\Delta t}}-\mu \frac{\partial ^{2}y^{n}}{\partial
x^{2}}+\epsilon \frac{\partial y^{n}}{\partial x}-y^{n}=\chi_{x_k}(x)
v_{k}^{n},$ \\
$-\mu \frac{\partial y^{n}\left( 0\right) }{\partial x}=v_{0}^{n},$ \\
$\mu \frac{\partial y^{n}\left( L\right) }{\partial x}=v_{M}^{n}.$%
\end{tabular}%
\ \right.  \tag{SE$^{\bigtriangleup t}$}  \label{SEdelT}
\end{equation}

The optimal condition is

\begin{equation*}
\delta J^{{\Delta t}}(v)=\sum_{k=0}^M\left( \nabla J^{^{{\Delta t}}}\left(
v_k\right) ,\delta v_k\right) _{\mathcal{U}^{{\Delta t}}}=0.
\end{equation*}

And
\begin{equation*}
\delta J^{{\Delta t}}(v)=k_0 \Delta t \sum_{k=0}^M
\sum_{n=0}^{N}v_k^{n}\delta v_k^{n}+k_{1}{\Delta t}\sum_{n=0}^{N}%
\int_{0}^{L}y^{n}\delta y^{n}\text{d}x+k_{2}\int_{0}^{L}y^{N+1}\delta y^{N+1}%
\text{d}x.
\end{equation*}

By the other hand, the perturbation of (\ref{SEdelT}) is

\begin{equation}
\left\{
\begin{tabular}{l}
$\delta y^{0}=0.$ \\
for $n=0,\ldots ,N$ \\
$\frac{\delta y^{n+1}-\delta y^{n}}{{\Delta t}}-\mu \frac{\partial
^{2}\delta y^{n}}{\partial x^{2}}+\epsilon \frac{\partial \delta y^{n}}{%
\partial x}-\delta y^{n}=\chi_{x_k}(x)\delta v_{k}^{n},$ \\
$-\mu \frac{\partial \delta y^{n}\left( 0\right) }{\partial x}=\delta
v_{0}^{n},$ \\
$\mu \frac{\partial \delta y^{n}\left( L\right) }{\partial x}=\delta
v_{M}^{n}.$%
\end{tabular}%
\ \right.  \tag{$\delta $SE$^{\Delta t}$}  \label{delSEdelT}
\end{equation}

Figure \ref{fig:delSE_DT} depicts (\ref{delSEdelT}).

Now, multiplying these by appropriate functions $p^{n}$ to integrate:
\begin{equation*}
{\Delta t}\sum_{n=0}^{N}\int_{0}^{L}p^{n}\left( \frac{\delta y^{n+1}-\delta
y^{n}}{{\Delta t}}-\mu \frac{\partial ^{2}\delta y^{n}}{\partial x^{2}}%
+\epsilon \frac{\partial \delta y^{n}}{\partial x}-\delta y^{n}-\chi
_{x_{k}}\delta v_{k}^{n}\right) \text{d}x=0.
\end{equation*}

\begin{equation}
{\Delta t}\sum_{n=0}^{N}\int_{0}^{L}p^{n}\left( \frac{\delta y^{n+1}-\delta
y^{n}}{{\Delta t}}\right) \text{d}x =  \label{6}
\end{equation}

\begin{eqnarray}
&=&-\int_{0}^{L}p^{0}\frac{\delta y^{0}}{{\Delta t}}\text{d}x-{\Delta t}%
\sum_{n=1}^{N}\int_{0}^{L}\left( \frac{p^{n}-p^{n-1}}{{\Delta t}}\right)
\delta y^{n}\text{d}x+\int_{0}^{L}p^{N}\delta y^{N+1}\text{d}x  \notag \\
&&-{\Delta t}\sum_{n=1}^{N}\int_{0}^{L}\left( \frac{p^{n}-p^{n-1}}{{\Delta t}%
}\right) \delta y^{n}\text{d}x+\int_{0}^{L}p^{N}\delta y^{N+1}\text{d}x.
\notag
\end{eqnarray}

\begin{equation}
{\Delta t}\sum_{n=0}^{N}\int_{0}^{L}p^{n}\left( -\mu \frac{\partial
^{2}\delta y^{n}}{\partial x^{2}}\right) \text{d}x = ~  \label{7}
\end{equation}

\begin{eqnarray}
&=&{\Delta t}\sum_{n=0}^{N}p\left[ -\mu \frac{\partial \delta y}{\partial x}%
\right] _{0}^{L}+\mu {\Delta t}\sum_{n=0}^{N}\frac{\partial p}{\partial x}%
\left[ \delta y\right] _{0}^{L}-\mu {\Delta t}\sum_{n=0}^{N}\int_{0}^{L}%
\frac{\partial ^{2}p}{\partial x^{2}}\delta y\text{d}x  \notag \\
&=&-{\Delta t}\sum_{n=0}^{N}p^{n}\left( 0\right) \left( \delta v^{n}\right)
+\mu {\Delta t}\sum_{n=0}^{N}\frac{\partial p^{n}\left( L\right) }{\partial x%
}\delta y\left( L\right)  \notag \\
&&-\mu {\Delta t}\sum_{n=0}^{N}\int_{0}^{L}\frac{\partial ^{2}p}{\partial
x^{2}}\delta y\text{d}x.  \notag
\end{eqnarray}

\begin{equation}
{\Delta t}\sum_{n=0}^{N}\int_{0}^{L}p^{n}\left( \epsilon \frac{\partial
\delta y^{n}}{\partial x}\right) \text{d}x=\epsilon
\sum_{n=0}^{N}p^{n}\left( L\right) \delta y\left( L\right) -\epsilon
\sum_{n=0}^{N}\int_{0}^{L}\frac{\partial p}{\partial x}\delta y\text{d}x.
\label{8}
\end{equation}

\begin{figure}
\centerline{ \psfig{figure=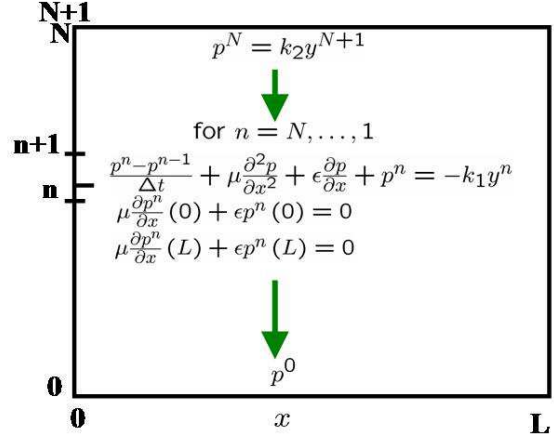,width=3.5in} } ~ ~
\caption{Discretization on time of adjoint system of (\protect\ref{SE}).}
\label{fig:SA_DT}
\end{figure}

\begin{equation}
{\Delta t}\sum_{n=0}^{N}\int_{0}^{L}p^{n}\left( -\delta y^{n}\right) \text{d}%
x.  \label{9}
\end{equation}

\begin{equation*}
0=(\ref{6})+(\ref{7})+(\ref{8})+(\ref{9})=
\end{equation*}

\begin{eqnarray*}
&& {\Delta t}\sum_{n=1}^{N}\int_{0}^{L}\left( -\frac{p^{n}-p^{n-1}}{{\Delta t%
}}-\mu \frac{\partial ^{2}p}{\partial x^{2}}-\epsilon \frac{\partial p}{%
\partial x}-p^{n}\right) \delta y^{n}\text{d}x+\int_{0}^{L}p^{N}\delta
y^{N+1}\text{d}x \\
&&-{\Delta t}\sum_{n=0}^{N}p\left( 0\right) \left( \delta v^{n}\right) +\mu {%
\Delta t}\sum_{n=0}^{N}\frac{\partial p^{n}}{\partial x}\left( L\right)
\delta y\left( L\right) +\epsilon \sum_{n=0}^{N}p\left( L\right) \delta
y\left( L\right) .
\end{eqnarray*}

Therefore the discretization on time of the adjoint system (see figure \ref%
{fig:SA_DT}) is

\begin{equation}
\left\{
\begin{tabular}{l}
$p^{N}=k_{2}y^{N+1}.$ \\
for $n=N,\ldots ,1$ \\
$\frac{p^{n}-p^{n-1}}{{\Delta t}}+\mu \frac{\partial ^{2}p}{\partial x^{2}}%
+\epsilon \frac{\partial p}{\partial x}+p^{n}=-k_{1}y^{n},$ \\
$\mu \frac{\partial p^{n}}{\partial x}\left( 0\right) +\epsilon p^{n}\left(
0\right) =0$ \\
$\mu \frac{\partial p^{n}}{\partial x}\left( L\right) +\epsilon p^{n}\left(
L\right) =0$.%
\end{tabular}%
\ \right.  \tag{ASE$^{\Delta t}$}  \label{ASEdelT}
\end{equation}

And
\begin{equation*}
\nabla J^{^{{\Delta t}}}\left( v\right) =\sum_{k=0}^{M}\left\{
v_{k}^{n}-p^{n}\left( 0\right) \right\} _{n=0}^{N}.
\end{equation*}

\subsection{Fully discretization}

Let $H>0,$ $H$ is an integer multiple of $M$, and $\bigtriangleup x={h}=%
\frac{L}{H}$. The indices for axis $x$ are $-1\leq j\leq H+1$. Note that two
sets of points are added on $j=-1$, and $j=H+1$, this is convenient because
the frontier conditions on $x=0$ ($-\mu \frac{\partial y\left( 0,t\right) }{%
\partial x}=v\left( t\right) ,$) and $x=L$ ($\mu \frac{\partial y\left(
L,t\right) }{\partial x}=0$) can be inserted before and after the points of
interest $0$ to $H$ on $x $.

\begin{figure}
\centerline{ \psfig{figure=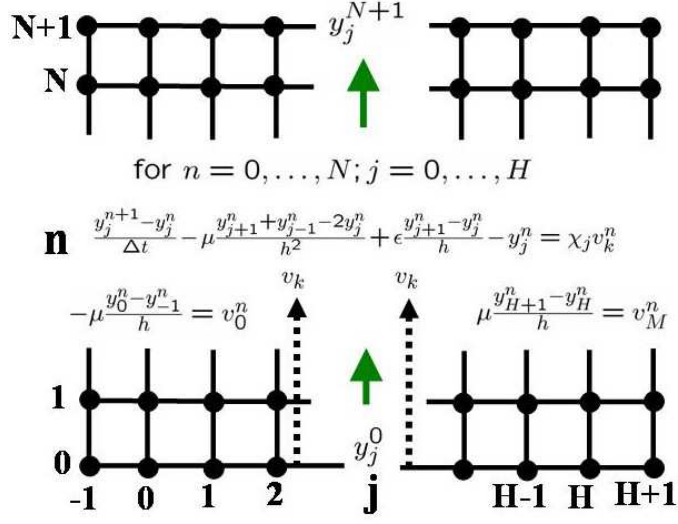,width=3.5in} } ~
\caption{Fully discretization of (\ref{SE}).} \label{fig:SE_DXDT}
\end{figure}

\begin{figure}
\centerline{ \psfig{figure=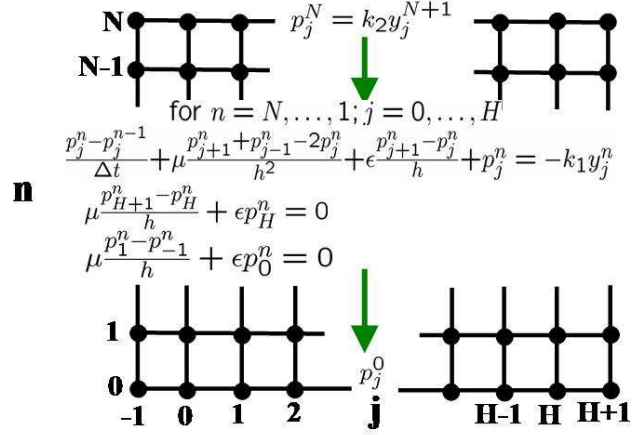,width=3.5in} } ~ ~
\caption{Fully discretization of adjoint system of (\ref{SE}).}
\label{fig:SA_DXDT}
\end{figure}

The corresponding fully discrete steady equations (see figure
\ref{fig:SE_DXDT})are

\begin{equation}
\left\{
\begin{tabular}{l}
$y_{j}^{0}=y_{0,j}$, $j=0,\ldots ,H$ \\
for $n=0,\ldots ,N$, $j=0,\ldots ,H$ \\
$\frac{y_{j}^{n+1}-y_{j}^{n}}{{\Delta t}}-\mu \frac{%
y_{j+1}^{n}+y_{j-1}^{n}-2y_{j}^{n}}{h^{2}}+\epsilon \frac{%
y_{j+1}^{n}-y_{j}^{n}}{h}-y_{j}^{n}=\chi _{x_{j}}v_{k}$ \\
$-\mu \frac{y_{0}^{n}-y_{-1}^{n}}{h}=v_{0}^{n}$ \\
$\mu \frac{y_{H+1}^{n}-y_{H}^{n}}{h}=v_{M}^{n}.$%
\end{tabular}%
\ \right.  \tag{SE$_{\vartriangle x}^{\vartriangle t}$}  \label{SEdelTdelX}
\end{equation}

$-\mu \frac{y_{0}^{n}-y_{-1}^{n}}{h}=v_{0}^{n}$

$-\mu y_{0}^{n}+\mu y_{-1}^{n}=hv_{0}^{n}$

$+y_{-1}^{n}=\left( hv_{0}^{n}+\mu y_{0}^{n}\right) /\mu $

$+y_{-1}^{n}=\frac{h}{\mu }v_{0}^{n}+y_{0}^{n}$

$\mu \frac{y_{H+1}^{n}-y_{H}^{n}}{h}=v_{M}^{n}$

$\mu \left( y_{H+1}^{n}-y_{H}^{n}\right) =hv_{M}^{n}$

$\mu y_{H+1}^{n}-\mu y_{H}^{n}=hv_{M}^{n} $

$\mu y_{H+1}^{n}=\mu y_{H}^{n}+hv_{M}^{n}$

$y_{H+1}^{n}=y_{H}^{n}+\frac{h}{\mu }v_{M}^{n}$

The adjoint equations (see figure \ref{fig:SA_DXDT}) are

\begin{equation}
\left\{
\begin{tabular}{l}
$p_{j}^{N}=k_{2}y_{j}^{N+1}$, $j=0,\ldots ,H$ \\
for $n=N,\ldots ,1$, $j=0,\ldots ,H$ \\
$\frac{p_{j}^{n}-p_{j}^{n-1}}{{\Delta t}}+\mu \frac{%
p_{j+1}^{n}+p_{j-1}^{n}-2p_{j}^{n}}{h^{2}}+\epsilon \frac{%
p_{j+1}^{n}-p_{j}^{n}}{h}+p_{j}^{n}=-k_{1}y_{j}^{n}$ \\
$\mu \frac{p_{H+1}^{n}-p_{H}^{n}}{h}+\epsilon p_{H}^{n}=0.$ \\
$\mu \frac{p_{0}^{n}-p_{-1}^{n}}{h}+\epsilon p_{-1}^{n}=0$%
\end{tabular}
\right.  \tag{ASE$_{\bigtriangleup x}^{\bigtriangleup t}$}
\label{ASEdelTdelX}
\end{equation}

$\frac{p_{j}^{n}-p_{j}^{n-1}}{{\Delta t}}+\mu \frac{%
p_{j+1}^{n}+p_{j-1}^{n}-2p_{j}^{n}}{h^{2}}+\epsilon \frac{%
p_{j+1}^{n}-p_{j}^{n}}{h}+p_{j}^{n}=-k_{1}y_{j}^{n}$.

The solution is $p_{j}^{n-1}=\frac{h^{2}p_{j}^{n}+\mu \Delta
tp_{j+1}^{n}+\mu \Delta tp_{j-1}^{n}-2\mu \Delta tp_{j}^{n}+\epsilon \Delta
thp_{j+1}^{n}-\epsilon \Delta thp_{j}^{n}+p_{j}^{n}\Delta
th^{2}+k_{1}y_{j}^{n}\Delta th^{2}}{h^{2}}=$

$p_{j}^{n-1}=p_{j}^{n}+\frac{\mu \Delta tp_{j+1}^{n}+\mu \Delta
tp_{j-1}^{n}-2\mu \Delta tp_{j}^{n}+\epsilon \Delta thp_{j+1}^{n}-\epsilon
\Delta thp_{j}^{n}+p_{j}^{n}\Delta th^{2}}{h^{2}}+k_{1}y_{j}^{n}\Delta t=$

$p_{j}^{n-1}=p_{j}^{n}+\frac{\mu \Delta t\left(
p_{j+1}^{n}+p_{j-1}^{n}-2p_{j}^{n}\right) }{h^{2}}+\frac{\epsilon \Delta
t\left( p_{j+1}^{n}-hp_{j}^{n}\right) }{h}+p_{j}^{n}\Delta
t+k_{1}y_{j}^{n}\Delta t=$

$\mu \frac{p_{0}^{n}-p_{-1}^{n}}{h}+\epsilon p_{-1}^{n}=0$

$\frac{\mu }{h}p_{0}^{n}-\frac{\mu }{h}p_{-1}^{n}+\epsilon p_{-1}^{n}=0$

$\mu p_{0}^{n}-\mu p_{-1}^{n}+\epsilon hp_{-1}^{n}=0$

$p_{-1}^{n}=\mu p_{0}^{n}/\left( \mu -\epsilon h\right) $

$\mu \frac{p_{H+1}^{n}-p_{H}^{n}}{h}+\epsilon p_{H}^{n}=0$

$\mu p_{H+1}^{n}-\mu p_{H}^{n}=-\epsilon hp_{H}^{n}$

$\mu p_{H+1}^{n}=\left( \mu -\epsilon \right) p_{H}^{n}/\mu. $

And the corresponding perturbation equations are

\begin{equation}
\left\{
\begin{tabular}{l}
$\delta y_{j}^{0}=0$, $j=0,\ldots ,H$ \\
for $n=0,\ldots ,N$, $j=0,\ldots ,H$ \\
$\frac{\delta y_{j}^{n+1}-\delta y_{j}^{n}}{{\Delta t}}-\mu \frac{\delta
y_{j+1}^{n}+\delta y_{j-1}^{n}-2\delta y_{j}^{n}}{h^{2}}+\epsilon \frac{%
\delta y_{j+1}^{n}-\delta y_{j}^{n}}{h}-\delta y_{j}^{n}=\chi _{x_{j}}\delta
v_{k}^{n}$ \\
$-\mu \frac{\delta y_{0}^{n}-\delta y_{-1}^{n}}{h}=\delta v_{0}^{n},$ \\
$\mu \frac{\delta y_{H+1}^{n}-\delta y_{H}^{n}}{h}=\delta v_{M}^{n}.$%
\end{tabular}%
\ \right.  \tag{$\delta $SE$_{\bigtriangleup x}^{\bigtriangleup t}$}
\label{delSEdelTdelX}
\end{equation}

The corresponding variational control problem is

\begin{equation}
\left\{
\begin{tabular}{l}
Find $u^\ast=\left\{ u^{n}\right\} \in \mathcal{V}_{\bigtriangleup
x}^{\bigtriangleup t} (=\mathbb{R}^{N \times M}$) \\
$J_{\bigtriangleup x}^{\bigtriangleup t}\left( u^\ast \right) \leq
J_{\bigtriangleup x}^{\bigtriangleup t}\left( v\right) \text{, }\forall\, v
\in \mathcal{V}_{\bigtriangleup x}^{\bigtriangleup t}$%
\end{tabular}%
\ \right.  \tag{CP$_{\bigtriangleup x}^{\bigtriangleup t}$}
\label{CPdelTdelX}
\end{equation}

where

\begin{equation*}
J_{\bigtriangleup x}^{\bigtriangleup t}\left( v\right) =\frac{{\Delta t}}{2}%
\sum_{k=0}^{M}\sum_{n=0}^{N}\left[ v_{k}^{n}\right] ^{2}+\frac{k_{1}{\Delta
th}}{2}\sum_{n=0}^{N}\sum_{j=0}^{H}\left[ y_{j}^{n}\right] ^{2}+\frac{k_{2}h%
}{2}\sum_{j=0}^{H}\left[ y_{j}^{N+1}\right] ^{2},
\end{equation*}

and $y=\left\{ y_{j}^{n}\right\} _{-1\leq j\leq H+1}^{0\leq n\leq N+1}$ is
the solution of (\ref{SEdelTdelX}) with $v$. Note that $H \geq M$ and $H$
must be a multiple of $M$ in order to have $\chi _{x_{j}}\delta v_{k}^{n}=
\delta v_{k}^{n}, \forall k=0,\ldots,M.$

\section{The Conjugate Gradient Algorithm}~\label{sc:CGALg}


The CG algorithm for the fully discrete control problem (\ref{CPdelTdelX})
is:

\begin{enumerate}
\item Given $\varepsilon $ (the tolerance to stop the algorithm), $%
0<\varepsilon \ll 1$, and $\ \left\{ u^{n,0}\right\}=\mathbf{0} \in $ $%
\mathcal{V}_{\bigtriangleup x}^{\bigtriangleup t}.$

\item Solve the equation (\ref{SEdelTdelX}), and

with the solution $\left\{ y_{j}^{n,0}\right\} _{-1\leq j\leq H+1}^{0\leq
n\leq N+1}$ solve (\ref{ASEdelTdelX}) to get $\left\{ p_{j}^{n,0}\right\}
_{-1\leq j\leq H+1}^{0\leq n\leq N}.$

\item Compute $g^{0}=\left\{ u_{j_{k}}^{n,0}+p_{j_{k}}^{n,0}\right\} _{0\leq
j_{k}\leq H}^{0\leq n\leq N}$, and set $w^{0}=g^{0}$.

Now, we have $u^{m}$, $g^{m}$, and $w^{m}$.

\item If $\frac{\left( g^{m+1},g^{m+1}\right) _{\mathcal{V}}}{\left(
g^{0},g^{0}\right) _{\mathcal{V}}}<\epsilon ^{2}$ take $u^{m+1}$ as the
solution and stop.

\item \label{GCciclo} Compute $m=m+1$.

\item Solve the equation ($\delta $SE$_{\bigtriangleup x}^{\bigtriangleup t}$%
), and

with the solution $\overline{y}=\left\{ \delta y_{j}^{n,m}\right\} _{-1\leq
j\leq H+1}^{0\leq n\leq N+1}$ solve (\ref{ASEdelTdelX}) to get $\overline{p}%
=\left\{ p_{j}^{n,m}\right\} _{-1\leq j\leq H+1}^{0\leq n\leq N}.$

\item Compute $\overline{g}^{m}=\left\{ w_{j_{k}}^{n,m}+\overline{p}%
_{j_{k}}^{n,m}\right\} _{0\leq j_{k}\leq H}^{0\leq n\leq N}$, $\rho
^{m}=\left( g^{m},g^{m}\right) _{\mathcal{V}}$, $u^{m+1}=u^{m}-\rho
^{m}w^{m} $, and $g^{m+1}=g^{m}-\rho ^{m}\overline{g}^{m}.$

\item If $\frac{\left( g^{m+1},g^{m+1}\right) _{\mathcal{V}}}{\left(
g^{0},g^{0}\right) _{\mathcal{V}}}<\epsilon ^{2}$ take $u^{m+1}$ as the
solution and stop.

\item Compute $\gamma ^{m}=\frac{\left( g^{m+1},g^{m+1}\right) _{\mathcal{V}}%
}{\left( g^{m},g^{m}\right) _{\mathcal{V}}}$, and $w^{m+1}=g^{m+1}+\gamma
^{m}w^{m}$

\item Go to step \ref{GCciclo}.
\end{enumerate}

\begin{figure}
\centerline{ \psfig{figure=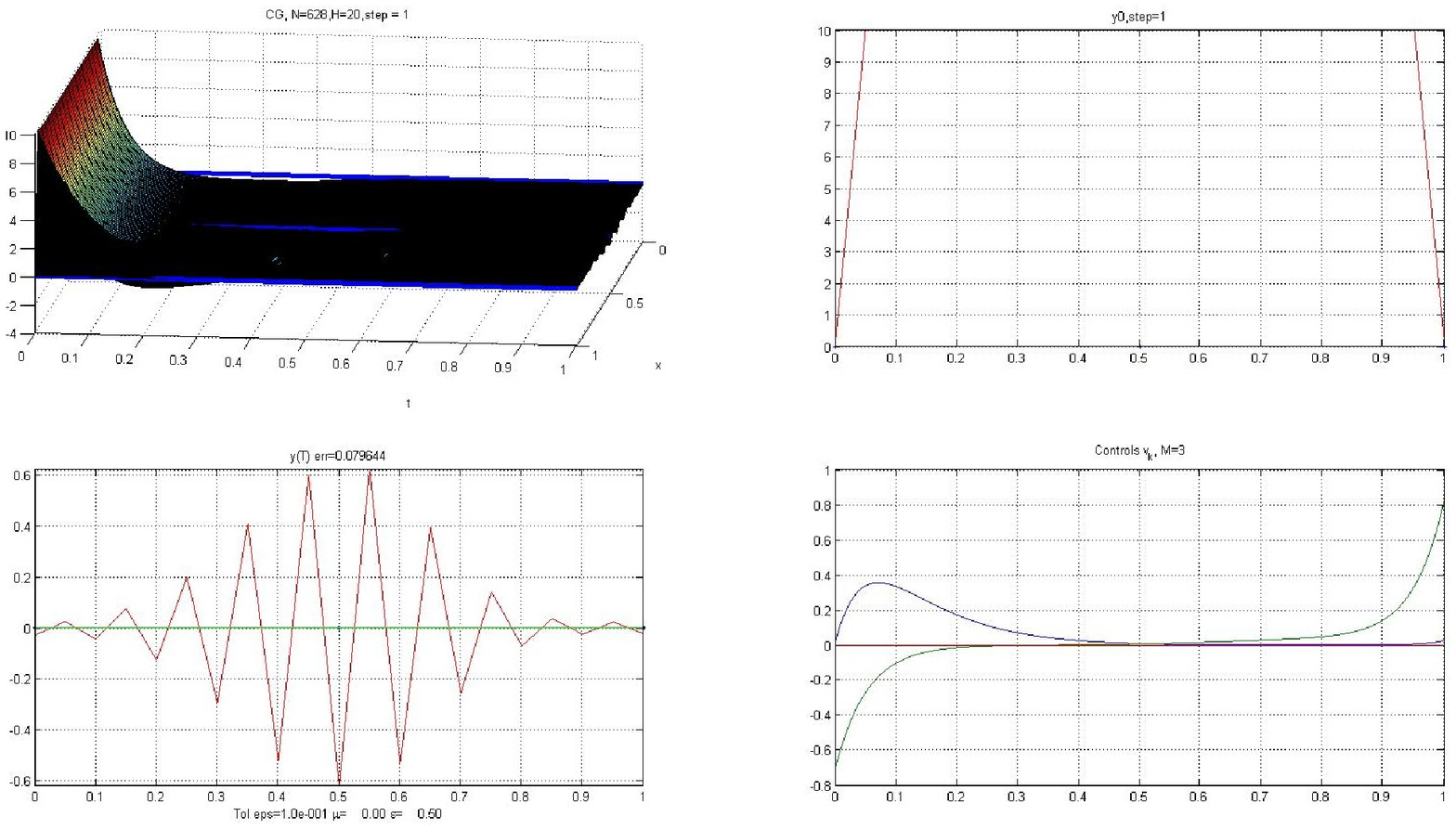,width=6.0in}}~ \caption{$y_0$
is a positive pulse, 3 controls.} \label{fig:P10M3}
\end{figure}

\begin{figure}
\centerline{ \psfig{figure=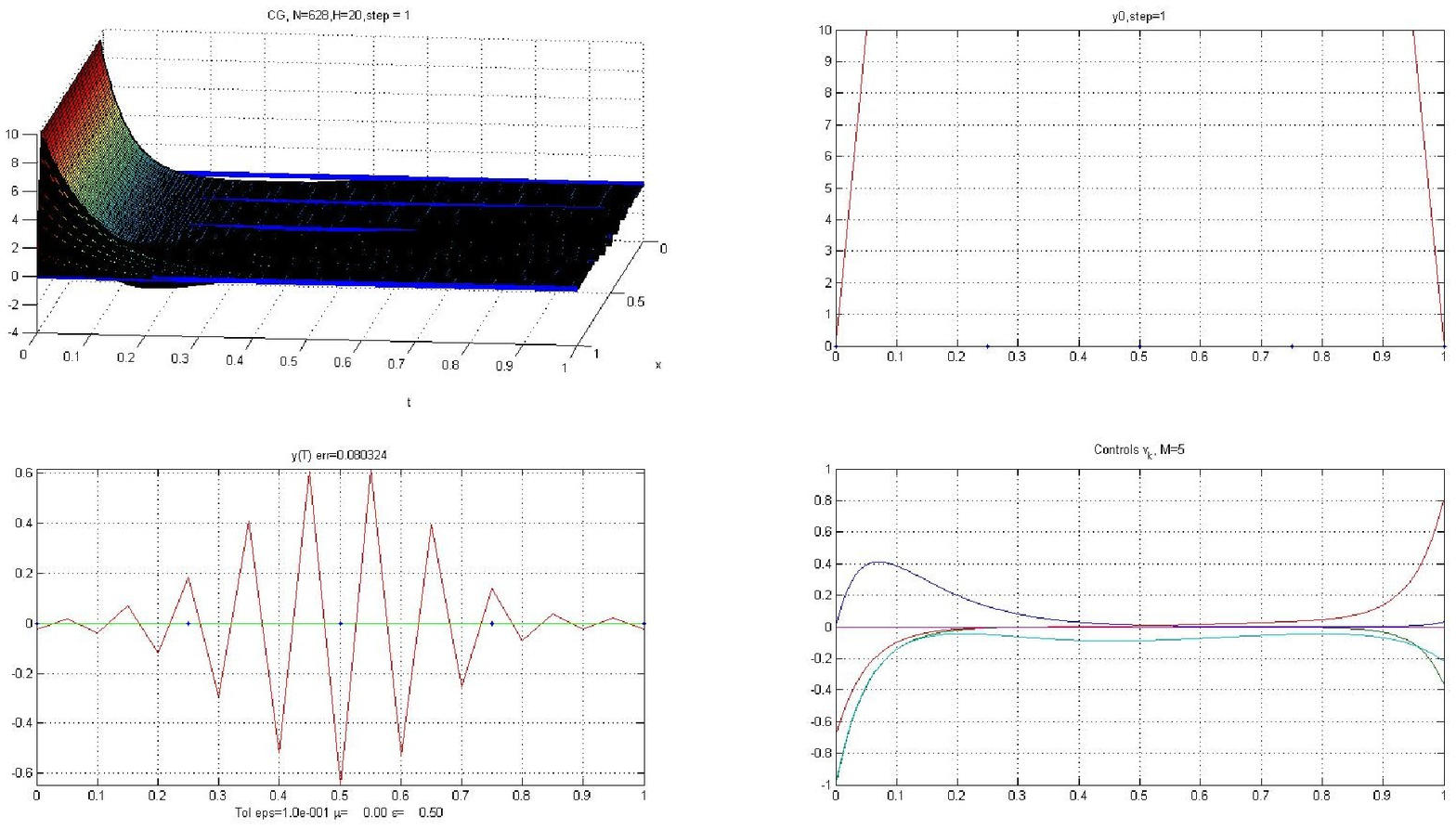,width=6.0in}}~ \caption{$y_0$
is a positive pulse, 5 controls.} \label{fig:P10M5}
\end{figure}

\begin{figure}
\centerline{ \psfig{figure=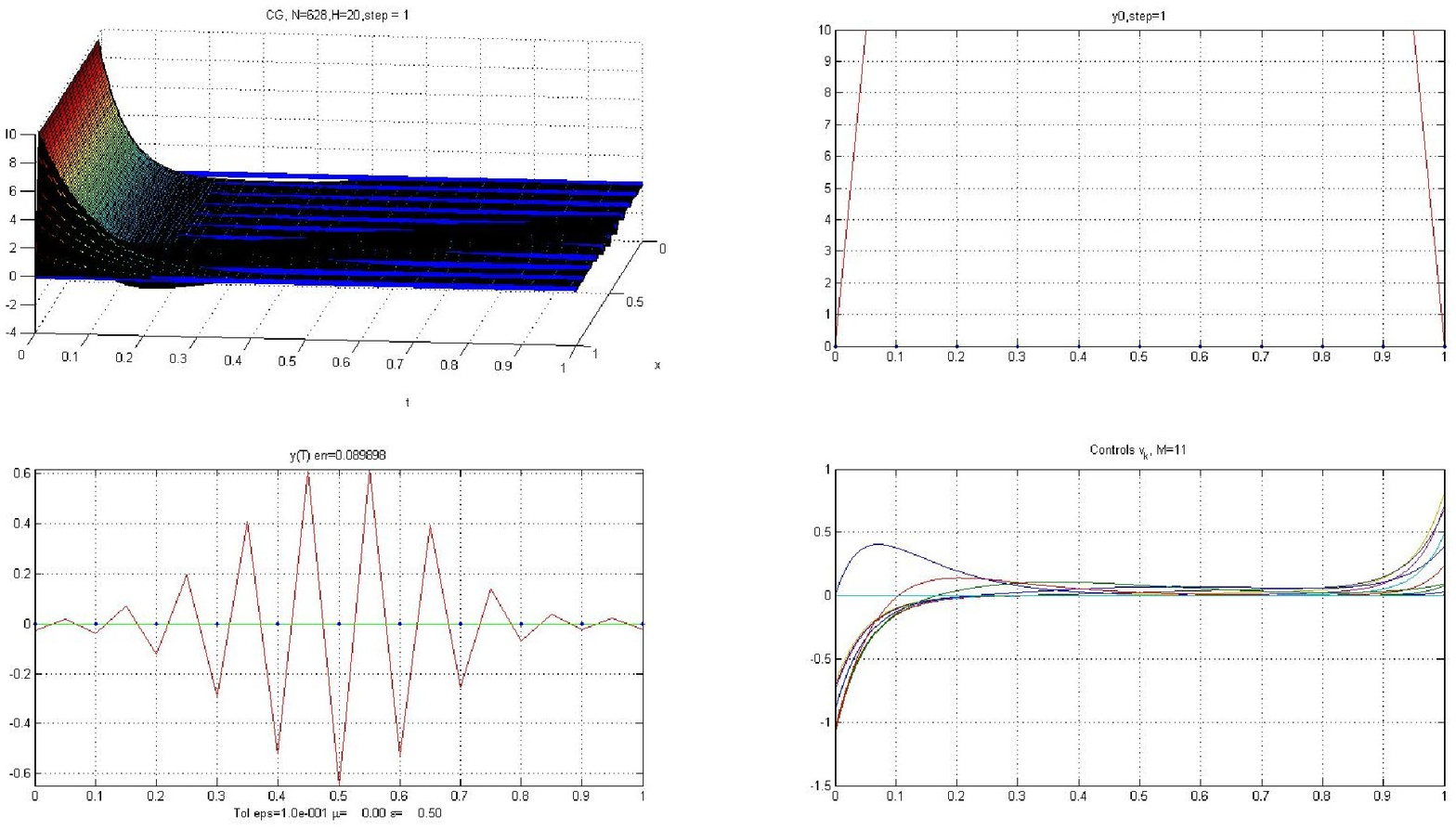,width=6.0in}}~ \caption{$y_0$
is a positive pulse, 11 controls.} \label{fig:P10M11}
\end{figure}

\section{Numerical Experiments}~\label{sc:numexp}

A program of the Conjugated Gradient Method
(Section~\ref{sc:CGALg}) was development in Matlab.

\begin{figure}
\centerline{ \psfig{figure=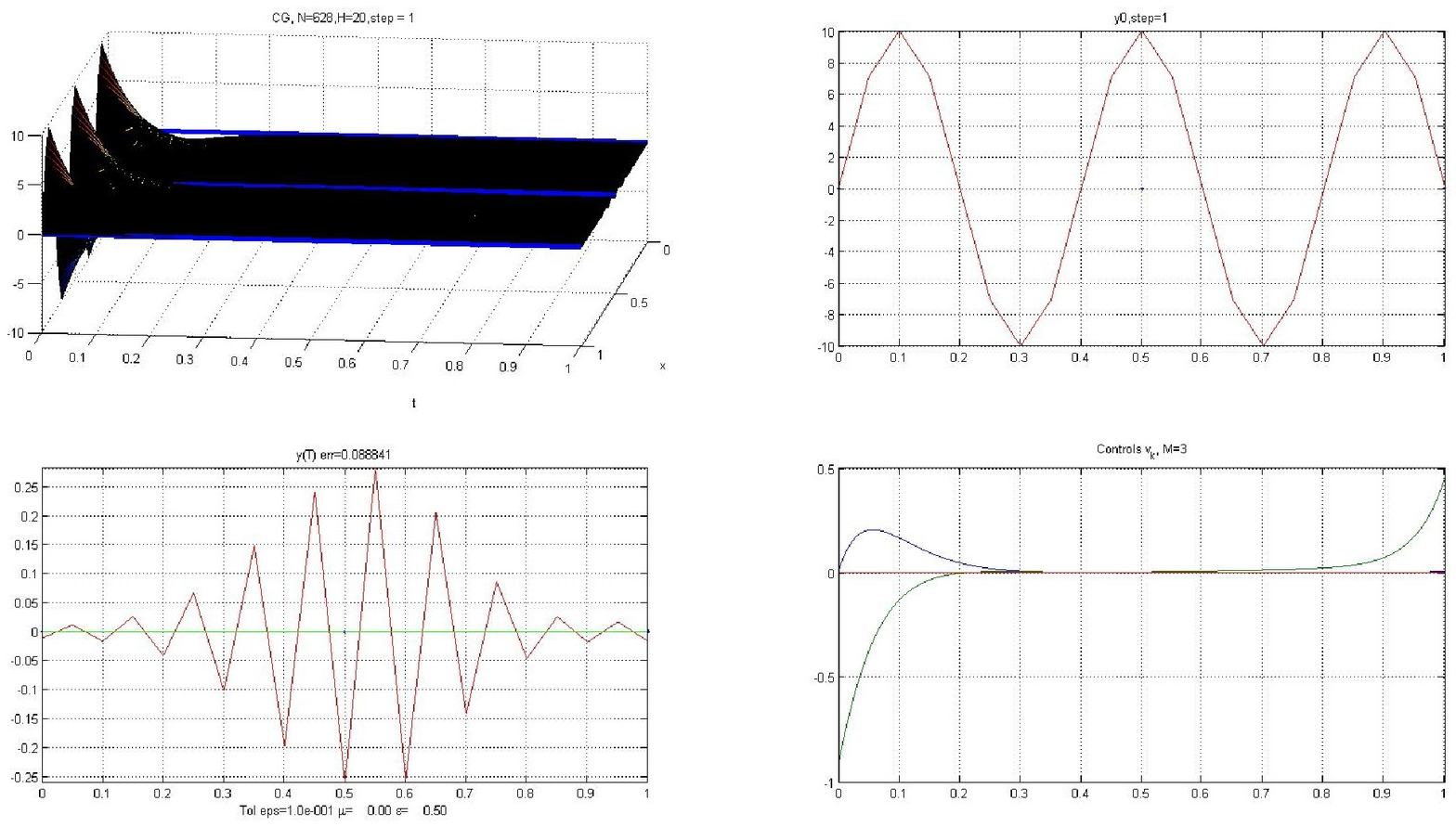,width=6.0in}}~
\caption{$y_0$ $=$ $10\sin(5\pi x),$ 3 controls.}
\label{fig:S2OA10M3}
\end{figure}

\begin{figure}
\centerline{ \psfig{figure=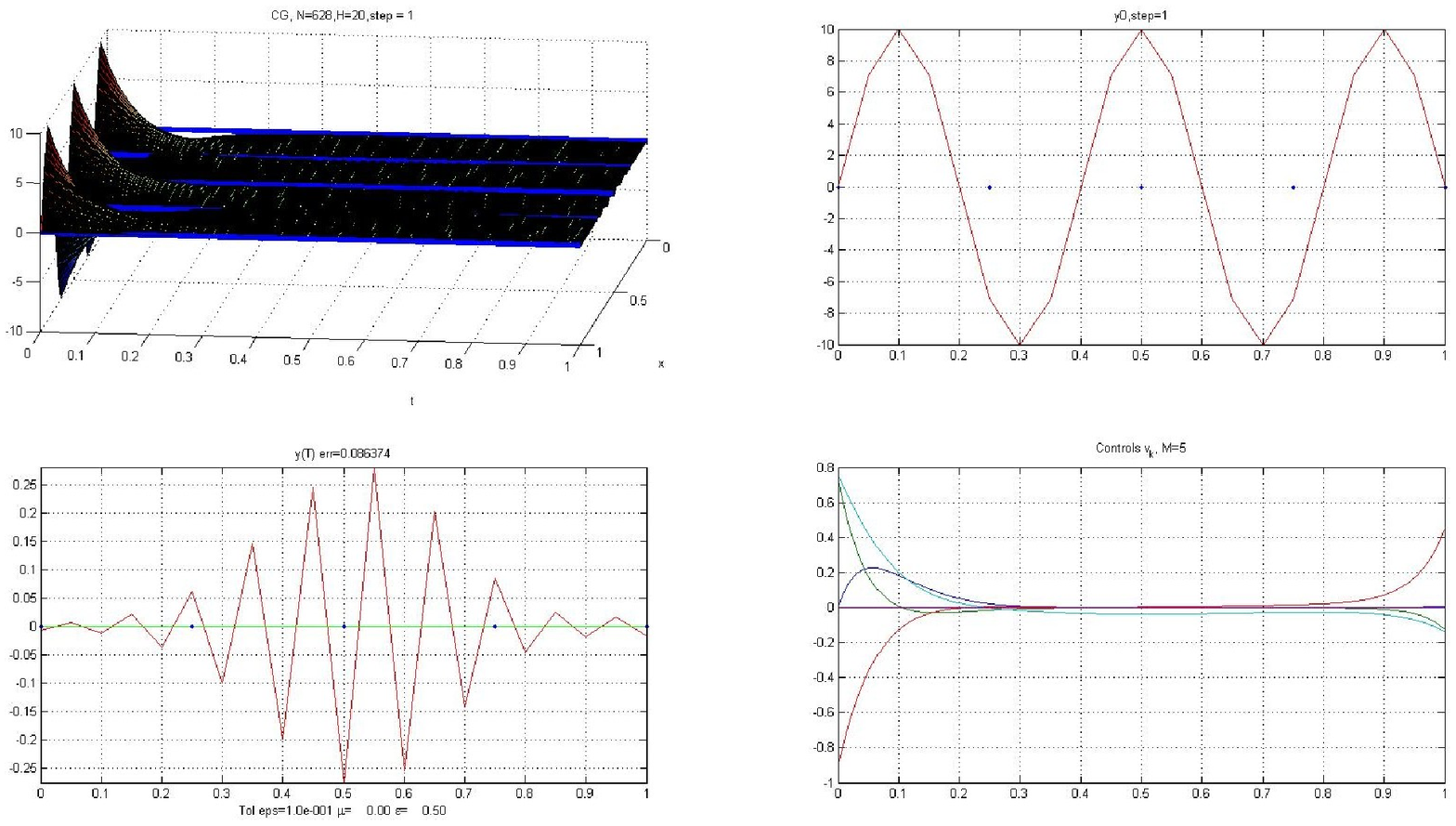,width=6.0in}}~
\caption{$y_0$ $=$ $10\sin(5\pi x),$ 5 controls.}
\label{fig:S2OA10M5}
\end{figure}

\begin{figure}
\centerline{ \psfig{figure=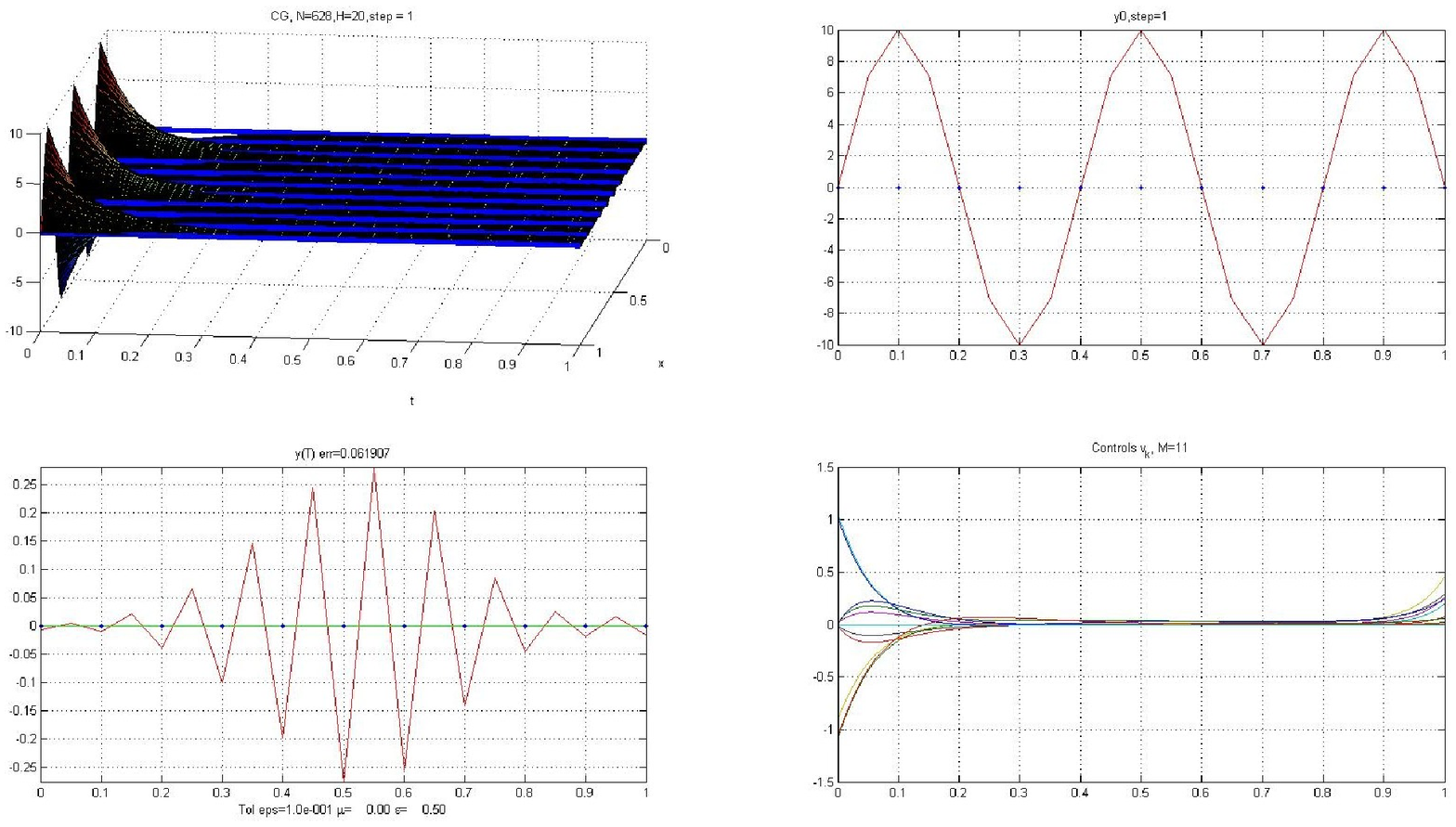,width=6.0in}}~
\caption{$y_0$ $=$ $10\sin(5\pi x),$ 11 controls.}
\label{fig:S2OA10M11}
\end{figure}

\begin{figure}
\centerline{ \psfig{figure=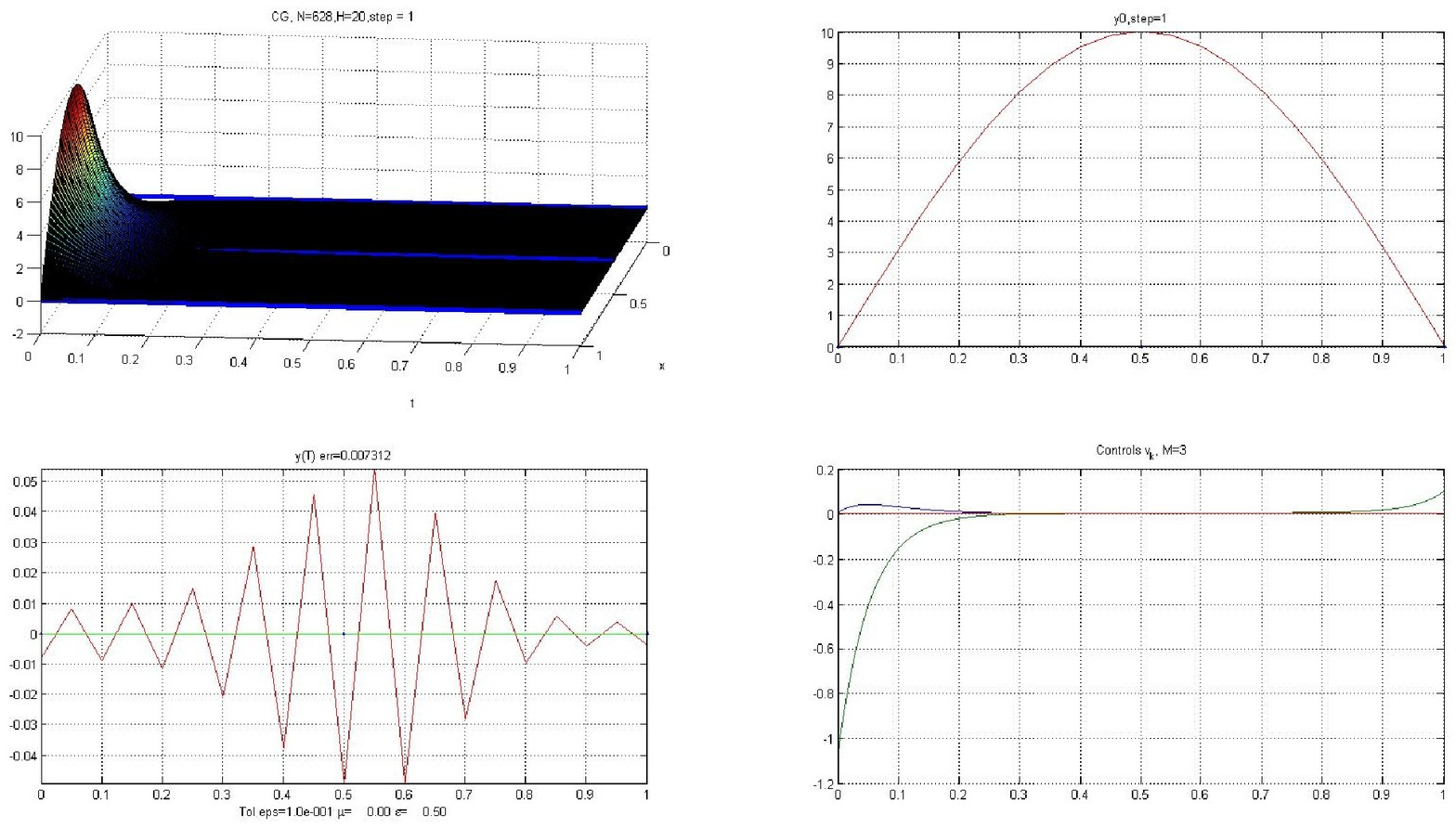,width=6.0in}}~ \caption{$y_0$
$=$ $10\sin(\pi x),$ 3 controls.} \label{fig:SA10M3}
\end{figure}

\begin{figure}
\centerline{ \psfig{figure=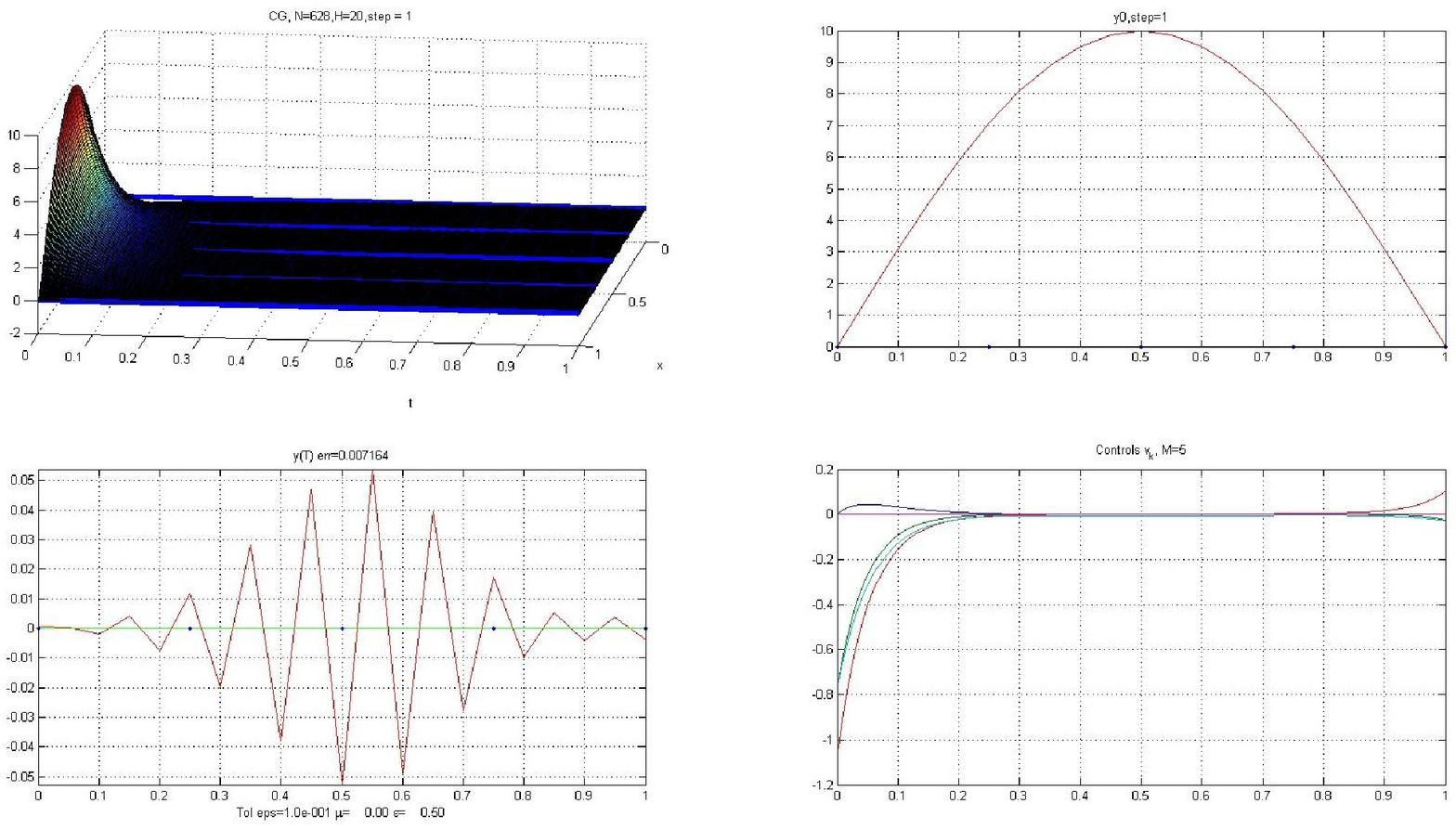,width=6.0in}}~ \caption{$y_0$
$=$ $10\sin(\pi x),$ 5 controls.} \label{fig:SA10M5}
\end{figure}

\begin{figure}
\centerline{ \psfig{figure=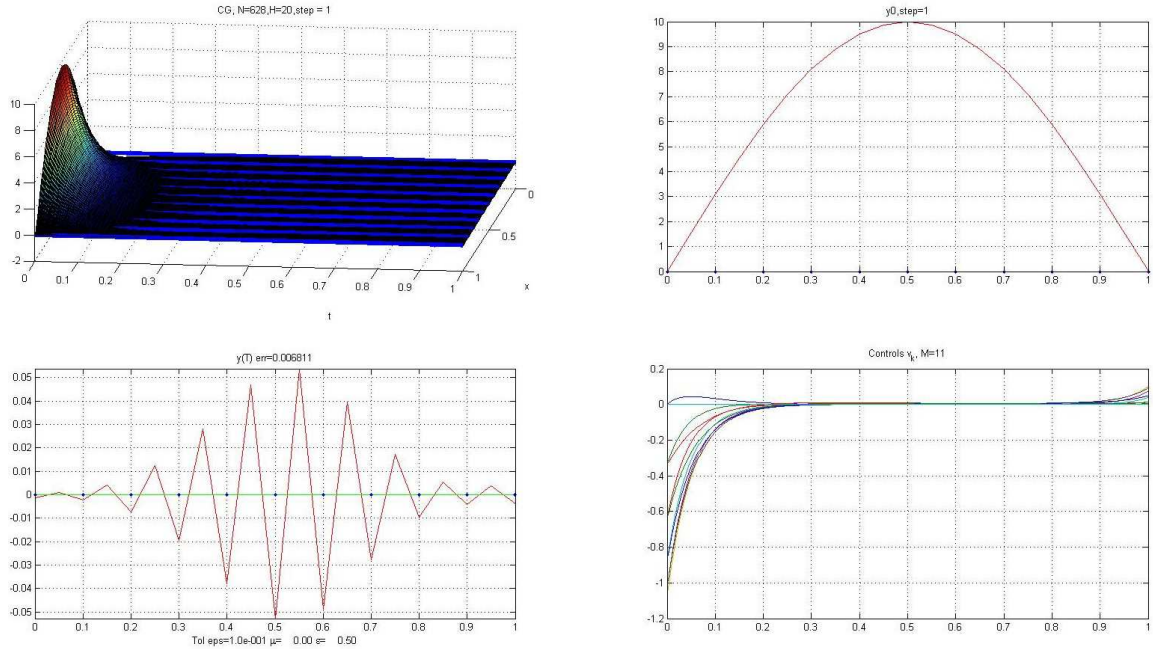,width=6.0in}}~
\caption{$y_0$ $=$ $10\sin(\pi x),$ 11 controls.}
\label{fig:SA10M11}
\end{figure}

Three numerical experiments were designed:
\begin{enumerate}
    \item $y_0$ is a positive pulse in $[0,1].$
    \item $y_0$ $=$ $10\sin(5\pi x),$ $x\in[0,1].$
    \item $y_0$ $=$ $10\sin(\pi x),$ $x\in[0,1].$

\end{enumerate}

The results are depicted in figures ~\ref{fig:P10M3},
~\ref{fig:P10M5}, ~\ref{fig:P10M11}, ~\ref{fig:S2OA10M3},
~\ref{fig:S2OA10M5}, ~\ref{fig:S2OA10M11}, ~\ref{fig:SA10M3},
~\ref{fig:SA10M5}, and ~\ref{fig:SA10M11}. Each figure depicts the
evolution of the state $y$, the initial state $y_0$, the final
state $y(T)$, and the graphs of the controls. These experiments
depict:

\begin{enumerate}
\item The graphs of the controls show that the cost increases with
the numbers of controls.

\item The graphs of the controls show that the controls behave
different.

\item It seems that with more controls the final state is closed
to steady state $\mathbf{0}$. However, there is a dependency of
the initial state $y_0$ and the numbers of controls in the
contrary. Figures~\ref{fig:P10M3}, ~\ref{fig:P10M5}, and
~\ref{fig:P10M11} depict a case where 3 controls behave better
than 11 controls.

\item It seems that with more controls the evolution of the $y$ is
controlled. In all cases, the controls are enough to diminish the
initial state $y_0$ and to keep under control the evolution of the
system over time.

\end{enumerate}


\section{Motivation for controlling}

We preferred to leave this section at the end, because these notes
are principally aimed for graduate students, which could be
interested in developing their own simulators. It is possibly,
that they already know the importance of the Theory of Control on
Systems over Partial Differential Equations or the control for
industrial process.

From the abundant literature, we mention the book of partial
differential
equations~\cite{cambrigePrees:Morton1994}, and for Control the books~\cite%
{springerCP:Glowinski1984, cambrigePrees:Glowinski2008}. These notes were
developed from the talk in~\cite{ CIMAT::Glowinski2006}.

The following problem depicts a classical problem for a parabolic equation
with three physical-chemical components.

\begin{enumerate}
\item Advection. It is the scalar variation at each point of a vector field,
by example, the contaminant entrainment in a medium.

\item Reaction. It is the response or reaction of the system, by example,
the heat exchanges in a system.

\item Diffusion. It is the gradient (change or transport) of system
components.
\end{enumerate}

\begin{figure}[ht]
\centerline{ \psfig{figure=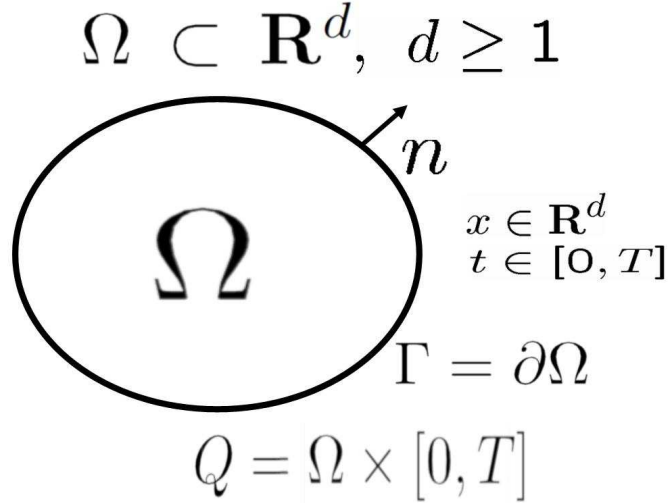,width=3.5in} } ~
\caption{Sistem's domain}
\label{fig:dominio}
\end{figure}

Let be the following parabolic equation where the advection is $V\cdot
\nabla \varphi $, the reaction is $f\left( \varphi \right) $, and the
diffusion is $\nabla \cdot \left( A\nabla \varphi \right) $ acting over the
time. It is Equation of the State System.

\begin{equation}
\begin{tabular}{l}
$\frac{\partial \varphi }{\partial t}-\nabla \cdot \left( A\nabla \varphi
\right) + V\cdot \nabla \varphi +f\left( \varphi \right) = 0 \text{ en }
Q=\Omega \times \left[ 0,T\right],$ \\
$A\nabla \varphi \cdot n = 0\text{ en }\Sigma =\Gamma \times \left[ 0,T%
\right],$ \\
$\varphi \left( x,0\right) = \varphi _{0}\left( x\right) \text{ }x\in \Omega
$%
\end{tabular}
~. ~  \tag{SEE}  \label{eq:SEE}
\end{equation}

where $\Omega \subset \mathbf{R}^{d}$ ($d\geq 1$, dimension) it is a smooth
region, with orientated boundary $\Gamma =\partial \Omega $, $n$ represents
a normal unit vector on $\Gamma $ (pointing outside of $\Omega $), $T>0$ is
the time ( including the possiblity $T=\infty )$. Figure~\ref{fig:dominio}
depicts~(\ref{eq:SEE}).

The intern product $\cdot $ is the usual, $a,b\in \mathbf{R}^{d},a\cdot
b=\sum\limits_{i=1}^{d}a_{i}b_{i}$, $A$ is a real tensor function (diffusion
matrix), $V:\Omega \rightarrow \mathbf{R}^{d}$ is a vectorial function, $f:%
\mathbf{R}\rightarrow \mathbf{R}$ is a real function, and $\varphi (x,t)$ is
the phenomena function that occurs in $Q.$

In addition we assume that:

\begin{equation*}
A(x)\xi \cdot \xi \geq \alpha \left\vert \xi \right\vert ^{2},\forall \xi
\in \mathbf{R}^{d}\text{ for almost all } x\in \Omega
\end{equation*}

which means that $A$ is uniformly positive definite for almost all $x $ in $%
\Omega.$

For the vector function $V,$ we assume:

\begin{eqnarray*}
\nabla \cdot V &=&0\text{ (divergence free)} \\
\frac{\partial V}{\partial t} &=&0\text{ (it is constant over time)} \\
V\cdot n &=&0\text{ on }\Gamma
\end{eqnarray*}

Control is necessary for this System, let be a reaction function given by

\begin{equation*}
f\left( \varphi \right) =C-\lambda e^{\varphi }
\end{equation*}

where $C,\lambda >0$ are real positive constants.

Then the steady state solution for such $f$ \ fulfill:

\begin{equation}
\frac{\partial \varphi }{\partial t} + f\left( \varphi \right) = 0
\label{eq:EPARDt}
\end{equation}

and it is given by

\begin{equation*}
\varphi _{s}=\frac{\ln C}{\lambda }
\end{equation*}

Note that $\varphi _{s}$ is constant,so that the equation~(\ref{eq:EPARDt}),
substituting $\varphi _{s}$ is fulfill (because $f\left( \varphi _{s}\right)
=C-\lambda e^{\varphi _{s}} = C-\lambda e^{\frac{\ln C}{\lambda }}=0)$.

Assuming that for some $t>0$, the system was its stable steady state
solution $\varphi =\varphi _{s}$.

Now, $\varphi =\varphi _{s}$ at some time $t_0=0$ has a small constant
perturbation $\delta \varphi ,$ independent from $x$ y $t$ (with $\nabla
\delta \varphi =0$ y $\frac{\partial \delta \varphi }{\partial t}=0$).

For this perturbation, the system evolves under the following ordinary
differential equation:

\begin{eqnarray*}
\frac{\text{d}\varphi }{\text{d}t} &=&\lambda e^{\varphi }-C\text{, }\lambda
,C>0, \text{ real constants} \\
\varphi \left( 0\right) &=&\varphi _{s}+\delta \varphi
\end{eqnarray*}

This model behaves with a constant positive perturbation, $\delta \varphi >0$%
, such that $\varphi \rightarrow +\infty.$ By other hand, if the
perturbation is a constant negative, $\delta \varphi <0$, then $\varphi
_{t\rightarrow \infty }\rightarrow -\infty $. In the following paragraphs,
it is showed that in the former the deviation from the stable state grows
fast to $+\infty$, and in the second case the deviation of the stable state
is slow and steady toward $-\infty$ as the time progress.

This means that around a stable steady state solution, the introduction a
small constant perturbation makes the system unstable. To verify the above
statement, we proceed by the Euler Method to numerically integrate the above
equation:

\begin{eqnarray*}
\frac{\text{d}\varphi }{\text{d}t} &=&\lambda e^{\varphi }-C\text{, }\lambda
,C>0,\text{ real constants} \\
\varphi \left( 0\right) &=&\frac{\ln C}{\lambda }+\delta \varphi
\end{eqnarray*}

Without loss of generality we take $\triangle t=1$, $C=1,\lambda =1,\delta
\varphi =0.1>0,$ and approach $\frac{d\varphi }{dt}$ by a time difference
between $n$ and $n-1$.

The resulting approximation difference equation is

$\varphi _{n}=\exp \left( \varphi _{n-1}\right) +\varphi _{n-1}-1.$

From the initial condition:

$\varphi _{0}=\frac{\ln C}{\lambda }+\delta \varphi =0.1$

The numerical estimations are

$\varphi _{1}=\exp \left( 0.1\right) +0.1-1=\allowbreak 0.205\,17$

$\varphi _{2}=\exp \left( 0.205\,17\right) +0.205\,17-1\allowbreak =$ $%
0.432\,9$

$\varphi _{3}=\exp \left( 0.432\,9\right) +0.432\,9-1=\allowbreak 0.974\,64$

$\varphi _{4}=\exp \left( 0.974\,64\right) +0.974\,64-1=\allowbreak
2.\,\allowbreak 624\,8$

$\varphi _{5}=\exp \left( 2.\,\allowbreak 624\,8\right) +2.\,\allowbreak
624\,8-1=\allowbreak 15.\,\allowbreak 427$

$\varphi _{6}=\exp \left( 15.\,\allowbreak 427\right) +15.\,\allowbreak
427-1=\allowbreak 5.\,\allowbreak 010\,3\times 10^{6}$

$\varphi _{7}=\exp \left( 5.\,\allowbreak 010\,3\times 10^{6}\right)
+5.\,\allowbreak 010\,3\times 10^{6}-1=\allowbreak 4.\,\allowbreak
392\,2\times 10^{2175945}$

$\varphi \left( t\right) $ in a finite time grows very quickly, it tends
accelerated to $\infty.$

By other hand, assuming that $\delta \varphi =-0.1<0$, and using the same
constants $C$ y $\lambda $, the numerical estimations for this case are

$\varphi _{0}=-0.1$

$\varphi _{1}=\exp \left( -0.1\right) +\left( -0.1\right) -1=\allowbreak
-0.195\,16$

$\varphi _{2}=\exp \left( -0.195\,16\right) +\left( -0.195\,16\right)
-1=\allowbreak -0.372\,46$

$\varphi _{3}=\exp \left( -0.372\,46\right) +\left( -0.372\,46\right)
-1=\allowbreak -0.683\,42$

$\varphi _{4}=\exp \left( -0.683\,42\right) +\left( -0.683\,42\right)
-1=\allowbreak -1.\,\allowbreak 178\,5$

$\varphi _{5}=\exp \left( -1.\,\allowbreak 178\,5\right) +\left(
-1.\,\allowbreak 178\,5\right) -1=\allowbreak -1.\,\allowbreak 870\,8$

$\varphi _{6}=\exp \left( -1.\,\allowbreak 870\,8\right) +\left(
-1.\,\allowbreak 870\,8\right) -1=\allowbreak -2.\,\allowbreak 716\,8$

$\varphi _{7}=\exp \left( -2.\,\allowbreak 716\,8\right) +\left(
-2.\,\allowbreak 716\,8\right) -1=\allowbreak -3.\,\allowbreak 650\,7$

$\varphi _{8}=\exp \left( -3.\,\allowbreak 650\,7\right) +\left(
-3.\,\allowbreak 650\,7\right) -1=\allowbreak -4.\,\allowbreak 624\,7$

$\varphi _{9}=\exp \left( -4.\,\allowbreak 624\,7\right) +\left(
-4.\,\allowbreak 624\,7\right) -1=\allowbreak -5.\,\allowbreak 614\,9$

$\varphi _{10}=\exp \left( -5.\,\allowbreak 614\,9\right) +-5.\,\allowbreak
614\,9-1=\allowbreak -6.\,\allowbreak 611\,3$

$\varphi _{11}=\exp \left( -6.\,\allowbreak 611\,3\right) +-6.\,\allowbreak
611\,3-1=\allowbreak -7.\,\allowbreak 610\,0$

$\varphi \left( t\right) $ is decreasing slowly to $-\infty .$

The previous numerical results clearly depicts that a control is necessary
to prevent such behavior and to return the system to the steady state
solution $\varphi _{s}$.

\section*{Conclusions and future work}~\label{sc:conclusions}

I did not expect implying that more controls means best result.
The numerical results depict this but the positive pulse. However,
in all numerical experiments the controls push back the controlled
system~\ref{SE} to the steady solution $\mathbf{0}$. As in global
optimization, the objective functions and problems have a relation
or compromise within the solution and the method for solving them.
Here, there are different behaviors between initial state and
numbers of controls.

The position of controls could be interesting to study in the
future.

My students of the master program in Engineering Process help me
to obtain preliminary results in less than three months. They
study the relation between one control and the initial state, they
found examples where one control does not work. Also, they want to
known about how difficult could be to apply advance mathematics
and to development control process software. I already have the
one control version, so they did the preliminary experiments. I
promise them, that I will development the multiple controls
version. My teaching philosophy is to help people to understand
and to be free of myths. Of course it is difficult but, it is
better to development a toy simulator than to buy one.

I believe, that is a good practice to help the students and people
to understand and take advance mathematics and to build by
themselves software, as an open box.

\section*{Acknowledgement}

Thanks to Adrian L\'{o}pez Ya\~{n}ez, Delia Rivera Ugalde, and Jos\'{e} \'{A}ngel
Sol\'{\i}s Herrera.

Dedicated to the 43 Ayotzinapa's students.


\end{document}